\documentclass[12pt,a4]{article}
\usepackage{amsmath,amssymb,amsfonts,graphicx,authblk}
\begin{document}

  \newtheorem{theorem}{Theorem}[section]
  \newtheorem{proposition}[theorem]{Proposition}
  \newtheorem{lemma}[theorem]{Lemma}
  \newtheorem{cor}[theorem]{Corollary}
  \newtheorem{remark}[theorem]{Remark}
  \newtheorem{algorithm}[theorem]{Algorithm}
\newcommand{\rdots}{.\hspace{.1em}\raisebox{.8ex}
                    {.\hspace{.1em}\raisebox{.8ex}{.}}}
\newcommand{\sC}
           {\mbox{\rm C\hspace{-0.4em}\rule{0.05ex}{1.55ex}\hspace{0.4em}}}
\newcommand{\sN}
           {\mbox{\rm N\hspace{-0.71em}\rule{0.05ex}{1.55ex}\hspace{0.71em}}}
\newcommand{\sR}
           {\mbox{\rm R\hspace{-0.7em}\rule{0.08ex}{1.55ex}\hspace{0.7em}}}
\newcommand{\sT}
           {\mbox{\rm T\hspace{-0.5em}\rule{0.05ex}{1.55ex}\hspace{0.5em}}}

\renewcommand{\theequation}{\arabic{section}.\arabic{equation}}

\newcommand{\be}{\begin{equation}}
\newcommand{\ee}{\end{equation}}
\newcommand{\ba}{\begin{array}}
\newcommand{\ea}{\end{array}}
\newcommand{\C}{{\mathbb C}}
\newcommand{\cH}{{\cal H}}
\newcommand{\al}{{\alpha}}
\newcommand{\bt}{{\beta}}
\newcommand{\G}{{\Gamma}}
\newcommand{\g}{{\gamma}}
\newcommand{\dl}{{\delta}}
\newcommand{\Dl}{{\Delta}}
\newcommand{\co}{{\rm const}\,}
\newcommand{\comment}[1]{}
\newcommand{\ds}{\displaystyle}
\newcommand{\eop}{\quad\rule{7pt}{8pt}\vspace*{3mm}}
\newcommand{\eq}[2]{\begin{equation}\label{#1}#2\end{equation}}
\newcommand{\ha}{\frac{1}{2}}
\newcommand{\La}{\Lambda}
\newcommand{\la}{\lambda}
\newcommand{\vp}{\varphi}
\newcommand{\tl}{\tilde}
\newcommand{\ma}[2]{\left[\begin{array}{#1}#2\end{array}\right]}
\newcommand{\Th}{\Theta}
\newcommand{\om}{\omega}
\newcommand{\Om}{\Omega}
\newcommand{\sg}{\sigma}
\newcommand{\sq}{\sqrt}
\newcommand{\proof}{{\sl Proof}.\ }
\newcommand{\txt}[1]{\quad{\rm #1}\quad}
\newcommand{\wh}{\widehat}
\newcommand{\wt}{\widetilde}

\newcommand{\canc}[1]{\textcolor[rgb]{.7,.7,.7}{#1}}
\newcommand{\ins}[1]{\textcolor[rgb]{1,0,0}{#1}}

\title{Computation of quasiseparable representations of Green matrices}

\date{}
\author[1]{P. Boito} 
\author[2]{Y. Eidelman}
\affil[1]{Dipartimento di Matematica, Universit\`a di Pisa, Largo Bruno Pontecorvo, 5 - 56127 Pisa, Italy. Email: paola.boito@unipi.it }
\affil[2]{School of Mathematical Sciences, Raymond and Beverly
Sackler Faculty of Exact Sciences, Tel-Aviv University, Ramat-Aviv,
69978, Israel. Email: eideyu@tauex.tau.ac.il}

\maketitle
\thispagestyle{empty}
\abstract{The well-known Asplund theorem states that the inverse of a (possibly one-sided) band matrix $A$ is a Green matrix. In accordance with quasiseparable theory, such a matrix admits a quasiseparable
representation in its rank-structured part. Based on this idea, we derive algorithms that compute a quasiseparable representation of $A^{-1}$ with linear complexity.

Many  inversion
algorithms for band matrices exist in the literature. However, algorithms based on a computation of the rank structure performed theoretically via the Asplund theorem
appear for the first time in this paper. Numerical experiments confirm complexity estimates and offer insight into stability properties.}

\section{Introduction}
\setcounter{equation}{0}

Inversion of band matrices is a central topic in structured linear algebra and has been extensively investigated in the literature \cite{Capovani70,YI79,BF81,Romani86, B91,BLR90,  BiniMeini99, Rosza91, Garcia2012, KS13}.
It has also sparked a wide range of results and generalizations, together with an interest towards rank-structured matrices, such as quasi/semi-separable and their variants; see e.g. \cite{GKbook, SSbiblio, VandebrilBook, EGH1, CG03} and references therein. In addition, let us mention that the interplay of Hessenberg-quasiseparable structure and polynomial recurrence relations is explored in \cite{Olshevsky10}: in particular, the authors give an explicit generator-based characterization of inverses of so-called twisted Green matrices, a subclass of quasiseparable matrices of order one.

Band matrices arise often in applications, for instance in connection with operators having a local action. Since an $n\times n$ banded matrix is completely defined by $O(n)$ parameters (for instance, its entries), we can expect that its inverse, if it exists, can also be defined by $O(n)$ parameters, even though it is typically full. This suggests that the inversion of banded matrices can be carried out with linear complexity, which is indeed the case.

The tridiagonal case, possibly in combination with additional structure, has been the object of detailed study \cite{Meurant92, Noschese13}. Conversely, it is natural to ask for a characterization of matrices whose inverse is tridiagonal. For nonsingular matrices, this problem is well-understood; in particular, it is known that the property of having a tridiagonal inverse is equivalent to (1,1)-semiseparability, i.e., a set of rank-one constraints on certain submatrices
\cite{VandebrilBook, Rosza91}.  Characterization of (possibly singular or rectangular) matrices whose generalized inverse is tridiagonal or, more generally, banded has also been investigated in the literature \cite{Bapat07, Bevilacqua05, Fasino02}. A recent contribution is found in \cite{BuenoFurtado}, where the authors focus on singular matrices $A$ whose Moore-Penrose inverse is irreducible and tridiagonal. Their approach relies on bordering techniques and yields a necessary and sufficient condition based on  rank properties not only of $A$, but also of bases of the null spaces of $A$ and $A^H$.

Aside from inversion algorithms, a different line of research focuses on {\em a priori} bounds on the inverses of banded matrices. As mentioned above, the inverse of a banded matrix $A$ is generally dense; but, under rather general hypotheses, it turns out that the entries of $A^{-1}$ decrease exponentially in absolute value w.r.t. their distance from the main diagonal. In other words, one can find constants $K>0$ and $0<\lambda<1$, dependent on bandwidth and spectral properties, but independent of $n$, such that
$$
|[A^{-1}]_{ij}|\leq K \lambda^{|i-j|}.
$$
The seminal paper on this topic \cite{DMS84} dates from 1984 and relies on polynomial approximation of the inverse function. Several improvements have been subsequently proposed, including extensions to more general matrix functions and to other sparsity patterns; see, e.g.,
\cite{BenziGolub99, BenziRazouk07, Canuto14, Frommer18, BenziRinelli} and references therein.

A fundamental tool in the analysis of inverses of banded matrices is the well-known Asplund theorem, which characterizes such inverses in terms of rank of certain submatrices. In this paper we revisit the problem of inverting banded matrices through quasiseparable structure \cite{EGH1}, in the light of Asplund's theorem. Indeed, the rank requirements on the inverse matrix -- that is, the property of being a Green matrix -- can be formulated in terms of quasiseparable generators. Starting from this representation, we establish algorithms that compute the quasiseparable generators of the inverse of a banded matrix, both in the one-sided (lower-banded) and in the two-sided (upper- and lower-banded) case. This is done in two ways: via QR factorization (Theorems \ref{UF1} and \ref{ICF1T}) and via LU factorization (Theorems \ref{UF1l} and \ref{I1CF1L}).
To the best of our knowledge, this specific approach to the inversion of banded matrices has not been pursued before in the literature.


The purpose of the present work is twofold. As already mentioned, the main contribution is a new explicit quasiseparable description of the structure of the inverse of a banded matrix, which brings together theoretical analysis and algorithm design. However, this result is also the starting point for further work about a quasiseparable viewpoint on decay bounds for inverses of banded matrices. Indeed, if a suitable explicit quasiseparable description of the inverse matrix is available, it is natural to ask whether this parameterization can be used, instead of polynomial approximation techniques, to give a tight characterization of the decay behavior. We expect that the results of Theorems \ref{UF1} and \ref{UF1l} will allow us to make progress along this research direction.

The paper is organized as follows. Section \ref{sec:background} recalls the 
Asplund theorem and the quasiseparable representation of Green matrices. 
Section \ref{sec:auxiliary} presents a general description of the inversion 
algorithm, which includes the factorization of a lower band matrix as a product 
of a transform matrix and an upper triangular one. The first factor turns out
to be a lower Green-upper band matrix represented as a product of elementary 
transformation matrices. The final result is obtained via multiplication of the
transform matrix by an upper triangular one and, because of invariance of Green 
matrices w.r.t. multiplication by triangular matrices, we obtain the product 
easily. The QR approach to matrix inversion, leading to Theorem \ref{UF1}, is 
detailed in Section \ref{sec:QR}, whereas its LU-based counterpart is 
presented in Section \ref{sec:LU}. Section \ref{sec:test} is devoted to 
numerical tests.


\section{Band and Green matrices, the block partitions}\label{sec:background}
\setcounter{equation}{0}

Let $r, N$ be integers such that $N>r>0$.
An $N\times N$ scalar matrix $A=\{A_{ij}\}_{i,j=1}^N$ is called {\it a lower
band matrix of order $r$}  if $A_{ij}=0$ for $i-j>r$.

A matrix $B$ is called {\it a lower Green matrix of order $r$} if
\be\label{natu}
{\rm rank}B(k:N,1:k+r-1)\le r,\quad k=1,2,\dots,N-r.
\end{equation}

A matrix $D$ is called {\it an upper Green matrix of order $r$} if
\be\label{natuu}
{\rm rank}D(1:k+r-1,k:N)\le r,\quad k=1,2,\dots,N-r.
\end{equation}

It is well known that the class of invertible lower Green of order $r$ matrices
coincides with the class of inverses of invertible lower band
matrices of the same order $r$.

\begin{theorem}[The Asplund theorem]
An invertible matrix $A$  is a lower 
band matrix of order $r$ if
and only if its inverse $B=A^{-1}$ is a lower Green matrix of order $r$.
\end{theorem}

The rank conditions (\ref{natu}) imply that  the matrix $B$ admits a
quasiseparable representation of a special type. Such
representations of matrices were studied in \cite{EGH1}. Let us recall the basic
definition. Let $F=\{F_{ij}\}_{i,j=1}^K$ be a block matrix
with entries of sizes $m_i\times n_j$. Assume that the strictly lower
triangular part of $F$ admits the representation
\be\label{qrprl}
F_{ij}=p(i)a(i-1)\cdots a(j+1)q(j),\;1\le j<i\le K
\end{equation}
where $p(i)\;(i=1,\dots,K),\;q(j)\;(j=1,\dots,K-1),\;
a(k)\;(k=2,\dots,K-1)$ are matrices of (small) sizes $m_i\times r_{i-1},r_i\times n_j,
r_k\times r_{k-1}$, respectively.

To get the quasiseparable representations of Green and band matrices we use the
block form of matrices. We treat $N\times N$ scalar matrices as $(N-r+2)\times
(N-r+2)$ block ones. To define the sizes
of the corresponding blocks we use the parameters
\be\label{msh}
\begin{gathered}
m_0=0,\;m_1=m_2=\dots=m_{N-r}=1,\;m_{N-r+1}=r;\\
n_0=r;\;n_1=n_2=\dots=n_{N-r}=1,\;n_{N-r+1}=0.
\end{gathered}
\end{equation}
Note that indices here start from zero.
A lower band of order $r$ matrix can be seen as a block one
with $n_i\times m_j,\;i,j=0,1,\dots,N-r+1$ entries and therefore turns out to be block upper
triangular. We treat a lower Green of order $r$ matrix as a block one with
entries of sizes $m_i\times n_j,\;i,j=0,1,\dots,N-r+1$.

Relative to this partition, the conditions (\ref{natu}) have the form
\be\label{natus}
{\rm rank}B(k+1:N-r+1,0:k)\le r,\quad k=0,1,2,\dots,N-r.
\end{equation}
We use the superscript $\prime$ to denote a block representation of a matrix.
Applying Theorem 5.9 and formula (4.9) in \cite{EGH1} we obtain in the
block form
$$
B^{\prime}(i,j)=p(i)a^>_{ij}q(j),\quad 0\le j<i\le N-r+1
$$
with matrices $p(i),\;i=1,\dots,N-r+1$ of sizes $m_i\times r$,
$q(j),\;j=0,\dots,N-r$ of sizes $r\times n_j$ and $a(k),\;k=1,\dots,N-r$
of sizes $r\times r$. Setting $j=s-1$ and omitting  the zero row and
the $N-r+1$-th column, which are empty matrices, we get
\be\label{lrz}
B(i,s-1)=p(i)a^>_{i,s-1}q(s-1),\quad 1\le s\le i\le N-r+1.
\end{equation}
Using Corollary 5.2 in \cite[p.87]{EGH1} we obtain the relations
\be\label{ola}
B(k:N,k-1)=P_kq(k-1),\quad k=1,\dots,N-r+1
\end{equation}
with $P_{N-r+1}=p^{\prime}(N-r+1),\;
P_k=\left(\ba{c}p(k)\\P_{k+1}a(k)\ea\right),\;k=N-r,\dots,1$.

Setting
\begin{gather*}
\tilde p(1)=p(1)q(0),\;\tilde q(0)=I_r,\;\tilde a(1)=a(1)q(0),\\
\tilde a(N-r)=p(N-r+1)a(N-r),\;\tilde p(N-r+1)=I_r,\;
\tilde q(N-r)=p(N-r+1)q(N-r)
\end{gather*}
we obtain the representation (\ref{lrz}) with $\tilde p(1),\tilde q(0),
\tilde a(1),\tilde p(N-r+1),\tilde q(N-r),\tilde a(N-r)$ instead of
$p(1),q(0),a(1),p(N-r+1),q(N-r),a(N-r)$. Hence without loss of generality one
can assume that in (\ref{lrz}) we have $q(0)=p(N-r+1)=I_r$. We assume always
that $q(0)=I_r$.

Thus the part $j-i\le r$ of a lower Green of order $r$ matrix, i.e. the block
strictly lower triangular part, is completely defined by the parameters
$p(i),q(i),a(i),\;i=1,\dots,N-r,\;p(N-r+1)$.
For instance an $(r+4)\times(r+4)$ matrix $B$ has the form
\begin{gather*}
B=\\
\left(\ba{ccccc}p(1)&\ast&\ast&\ast&\ast\\
p(2)a(1)&p(2)q(1)&\ast&\ast&\ast\\
p(3)a(2)a(1)&p(3)a(2)q(1)&p(3)q(2)&\ast&\ast\\
p(4)\cdots a(1)&p(4)a(3)a(2)q(1)&p(4)a(3)q(2)&p(4)q(3)&\ast\\
p(5)\cdots a(1)&p(5)a(4)a(3)a(2)q(1)&p(5)a(4)a(3)q(2)&
p(5)a(4)q(3)&p(5)q(4)\ea\right).
\end{gather*}
The elements $p(i),q(i),a(i),\;i=1,\dots,N-r,\;p(N-r+1)$, where
$p(i)\;(i=1,\dots,N-r)$ are $r$-dimensional rows and $p(N-r+1)$ is an
$r\times r$ matrix, $q(i)$ are $r$-dimensional
columns and $a(i)$ are $r\times r$ matrices,
are said to be {\it lower Green generators} of the matrix $B$. One can check
easily (see \cite[Lemma 5.8]{EGH1}) that if the representation (\ref{lrz})
holds then $B$ is a lower Green of order $r$ matrix.

\section{The inversion of lower band matrices}\label{sec:auxiliary} 
\setcounter{equation}{0}

\subsection{The transform and the factorization}

For a lower band matrix of order $r$ 
\be\label{marf}
A=\left(\ba{ccccc}a_{11}&a_{12}&a_{13}&\dots&a_{1,N}\\
a_{21}&a_{22}&a_{23}&\dots&a_{2,N}\\
\vdots&\vdots&\vdots&\ddots&\vdots\\
a_{r+1,1}&a_{r+1,2}&a_{r+1,3}&\dots&a_{r+1,N}\\
0&a_{r+2,2}&a_{r+2,3}&\dots&a_{r+2,N}\\
0&0&a_{r+3,3}&\cdots&a_{r+3,N}\\
\vdots&\vdots&\vdots&\ddots&\vdots\\
0&0&0&\dots&a_{NN}\ea\right)
\end{equation}
we determine the $(r+1)\times(r+1)$ transforms $G_i,\;i=1,\dots,N-r+1$, which 
reduce $A$ to an upper triangular form, i.e., the $N\times N$ matrix $G$ such 
that 
$$
GA=R,
$$
with an upper triangular $R$, and next compute
$$
A^{-1}=R^{-1}G.
$$
The transform matrix $G$ is a product of elementary transform matrices. Using 
the lower band form (\ref{marf}) of the matrix $A$, the matrix $G$ may be written as the product 
\be\label{unb}
G=\tilde G_{N-r+1}\tilde G_{N-r}\tilde G_{N-r-1}\cdots\tilde G_1
\end{equation}
with
\be\label{unb1}
\tilde G_k=I_{k-1}\oplus G_k\oplus I_{N-k-r},\; k=1,\dots,N-r,\quad
\tilde G_{N-r+1}=I_{N-r}\oplus G_{N-r+1}.
\end{equation}
Note that the number of factors is linear in $N$.

It turns out that the matrix $G$ defined in (\ref{unb}), (\ref{unb1}) is lower 
Green and upper band  with the same order at the same time. Moreover, the
lower Green generators of $G$ may be obtained easily.

\begin{lemma}\label{UU}
Let $G$ be an $N\times N$ matrix which admits the factorization (\ref{unb}),
(\ref{unb1}),
where $G_k,\;k=1,\dots,N-r$ are $(r+1)\times(r+1)$ matrices and $G_{N-r+1}$
is an $r\times r$ matrix. Assume that the matrices
$G_k\;k=1,2,\dots N-r$ are partitioned in the form
\be\label{iof5b}
G_k=\left[\ba{cc}p_G(k)&d_G(k)\\a_G(k)&q_G(k)\end{array}\right],\quad
k=1,\dots,N-r,\quad G_{N-r+1}=p_G(N-r+1)
\end{equation}
with submatrices $p_G(k),d_G(k),a_G(k),q_G(k)$ of sizes $1\times r,
1\times 1,r\times r,r\times 1$ respectively, and $r\times r$ matrix
$p_G(N-r+1)$.

Then $G$ is a lower Green and upper band of order $r$ matrix with lower Green
generators $p_G(i),q_G(i),a_G(i),\;i=1,\dots,N-r,\;p_G(N-r+1)$ and
diagonal in the block form entries $G(k,k+r)=d_G(k),\;k=1,\dots,N-r$.
\end{lemma}

{\em Proof.} We apply Lemma 20.2 in \cite[p.375]{EGH1} to the matrix $G$. We set
$G_0=I_{r+1}$ and $\tilde G_0=I_N$ and, using (\ref{unb}), we obtain the 
representation
$$
G=\tilde G_{N-r+1}\tilde G_{N-r}\tilde G_0.
$$
We treat the matrix $G$ in a block form with entries of sizes 
$m_i\times n_j,\;i,j=0,1,\dots,N-r+1$ with $m_i,n_j$ as in (\ref{msh}).
Lemma 20.2 in \cite[p.375]{EGH1} implies that the block matrix $G$ is upper
triangular with lower quasiseparable generators and block diagonal entries
obtained from the partitions
$G_0=\left[\ba{c}d_G(0)\\q_G(0)\end{array}\right]$ and (\ref{iof5b}). The 
equality
$n_0=r$ implies that $G$ is an upper band matrix of order $r$. Moreover, using
(\ref{lrz}), we conclude that $p_G(i),q_G(i),a_G(i),\;i=1,\dots,N-r,\;
p_G(N-r+1)$ are lower Green
generators of $G$.$\hfill\Box$

The reverse statement is also true.

\begin{lemma}\label{UU1}
Let $G$ be an $N\times N$ lower Green and upper band of order $r$ matrix with 
lower Green generators $p_G(i),q_G(i),a_G(i),\;i=1,\dots,N-r,\;p_G(N-r+1)$ and
and the entries 
$$
d_G(k)=G(k,k+r),\;k=1,\dots,N-r.
$$
Set
\be\label{iof5b1}
G_k=\left[\ba{cc}p_G(k)&d_G(k)\\a_G(k)&q_G(k)\end{array}\right],\quad
k=1,\dots,N-r,\quad G_{N-r+1}=p_G(N-r+1).
\end{equation}

Then the matrix $G$ satisfies the formula (\ref{unb}) with the factors $\tilde
G_k$ as in (\ref{unb1}).
\end{lemma}

The proof follows directly from Lemma 20.1 in  \cite[p.374]{EGH1}.

{\bf Remark.} The formula (\ref{unb}) implies that the relation
\be\label{nati}
W=\tilde W_0\tilde W_1\cdots\tilde W_{N-r+1},
\end{equation}
with $\tilde W_i$ as  in (\ref{unb1}), i.e.
\be\label{unb1u}
\tilde W_k=I_{k-1}\oplus W_k\oplus I_{N-k-r},\; k=1,\dots,N-r,\quad
\tilde W_{N-r+1}=I_{N-r}\oplus W_{N-r+1}
\end{equation}
where $W_k,\;k=1,\dots,N-r$ are $(r+1)\times(r+1)$ matrices and $W_{N-r+1}$
is an $r\times r$ matrix, yields an upper Green and lower band of order $r$
matrix.  As an  example (with $r=1$) one can take
a unitary Hessenberg matrix. Such representations for the case $r=1$ for 
Green matrices as 
well as the permutations of (\ref{nati}) leading to the CMV and to the Fiedler
matrices, have been studied by V.~Olshevsky, G.~Strang and P.~Zhlobich in the paper
\cite{Olshevsky10}.

\subsection{The multiplication by triangular matrix}
Next we compute the inverse matrix $A^{-1}$ via the formula $A^{-1}=R^{-1}G$.  
To this end we derive an algorithm of multiplication of a Green matrix by an
upper triangular one.

\begin{theorem}\label{ULM}
Let $S$ be an upper triangular matrix and $B$ be a lower Green of order $r$
matrix with lower Green generators $p_B(i),q_B(i),
a_B(i)$, $i=1,\dots,N-r$, $p_B(N-r+1)$.

Then $C=SB$ is a lower Green of order $r$ matrix with lower Green generators
$$
p_C(i),q_C(i),a_C(i),\;i=1,\dots,N-r,\quad p_C(N-r+1)
$$
obtained as follows.

1. Set
\be\label{mas}
q_C(k)=q_B(k),\;a_C(k)=a_B(k),\;k=1,\dots,N-r
\end{equation}
and
\be\label{lrs}
\begin{gathered}
s_k=S(k,k),\;S_k=S(k:N,k:N),\;k=1,\dots,N,\\ l_k=S(k,k+1:N),\;k=1,\dots,N-1.
\end{gathered}
\end{equation}

2. Set
\be\label{mishn}
P^B_{N-r+1}=p_B(N-r+1),\;p_C(N-r+1)=S_{N-r+1}p_B(N-r+1),
\end{equation}
and for $k=N-r,\dots,1$ compute
\be\label{mish}
p_C(k)=s_kp_B(k)+l_kP^B_{k+1}a_B(k).
\end{equation}
\be\label{mash}
P^B_k=\left(\ba{c}p_B(k)\\P^B_{k+1}a_B(k)\ea\right).
\end{equation}
\end{theorem}

{\em Proof.}
Using (\ref{ola}) we get
\be\label{lenns}
B(k:N,k-1)=P^B_kq_B(k-1),\quad k=1,\dots,N-r+1
\end{equation}
with  $P^B_k$ as in (\ref{mishn}), (\ref{mash}).
Hence we have
\begin{gather*}
C(k:N,k)=S(k:N,:)B(:,k)=\\
\left(\ba{cc}S(k:N,1:k-1)&S_k\ea\right)
\left(\ba{c}B(1:k-1,k)\\B(k:N,k)\ea\right),\\ k=1,\dots,N-r+1
\end{gather*}
and using (\ref{lenns}) and the fact that $S$ is an upper triangular matrix we
get
$$
C(k:N,k)=S_kP^B_kq_B(k-1),\quad k=1,\dots,N-r+1.
$$
We have obviously
$$
S_k=\left(\ba{cc}s_k&l_k\\0&S_{k+1}\ea\right),\;k=1,\dots,N-r.
$$
Set also $P^C_k=S_kP^B_k,\;k=1,\dots,N-r+1$. We have
$$
P^C_{N-r+1}=S_{N-r+1}P^B_{N-r+1}=S_{N-r+1}p_B(N-r+1)=p_C(N-r+1).
$$
and next
$$
P^C_k=\left(\ba{cc}s_k&l_k\\0&S_{k+1}\ea\right)
\left(\ba{c}p_C(k)\\P^C_{k+1}a_B(k)\ea\right)
$$
with $p_C(k)$ as in (\ref{mish}). Thus we obtain that
$$
C(k:N,k)=P^C_kq_C(k-1),\quad k=1,\dots,N-r,
$$
where
$$
P^C_{N-r+1}=p_C(N-r+1),\quad
P^C_k=\left(\ba{c}p_C(k)\\P^C_{k+1}a_C(k)\ea\right),
\;k=N-r-1,\dots,1
$$
with $p_C(k),q_C(k),a_C(k)$ as in (\ref{mish}) and (\ref{mas}). Hence using
Lemma 5.3  in  \cite[p.88]{EGH1} we conclude that $p_C(k),q_C(k),a_C(k),\;
k=1,\dots,N-r,\;p_C(N-r+1)$ are lower Green generators of the matrix $C$. 
$\hfill\Box$

\section{The QR algorithm}\label{sec:QR}
\setcounter{equation}{0}

We consider here the case where $X$ is a unitary matrix. The main contributions of this section are Theorem \ref{UF1}, which describes the inversion algorithm for lower-banded matrices, and Theorem \ref{ICF1T}, which concerns inversion of two-sided banded matrices.

\subsection{Inverses of lower band matrices}

Let us begin with a characterization of the unitary-triangular factorization of a lower banded matrix.

\begin{theorem}\label{ICF1}
Let $A=\{A_{ij}\}_{i,j=1}^N$ be a lower band of order $r$ matrix.

The matrix $A$ admits the factorization
\be\label{irn15}
A=UR,
\end{equation}
where $U$ is a unitary matrix represented as the product
\be\label{un}
U=\tilde U_1\tilde U_2\cdots\tilde U_{N-r}\tilde U_{N-r+1}
\end{equation}
with
\be\label{un1}
\tilde U_k=I_{k-1}\oplus U_k\oplus I_{N-k-r},\; k=1,\dots,N-r,\quad
\tilde U_{N-r}=I_{N-r}\oplus U_{N-r+1},
\end{equation}
where $U_k,\;k=1,\dots,N-r$ are $(r+1)\times(r+1)$ unitary matrices,
is an $r\times r$ unitary matrix  , and $R$ is an upper triangular
matrix. Moreover the unitary matrices $U_k$ as well as the upper triangular
entries of the matrix $R$ are obtained as follows.

1. Set
\be\label{lry0}
Y_0=A(1:r,1:N).
\end{equation}

2. For $k=1,\dots,N-r$ perform the following.

Set
\be\label{lnk}
\Delta_k=\left[\ba{c}Y_{k-1}(:,1)\\A(k+r,k)\ea\right]
\end{equation}
and determine an $(r+1)\times(r+1)$ unitary matrix $U_k$ and a number $x_k$
such that
\be\label{lenk}
U_k^*\Delta_k=\left[\ba{c}x_k\\0_{r\times1}\ea\right].
\end{equation}
Compute the $(r+1)\times(N-k)$ matrix
\be\label{ol}
Z_k=U_k^*\left[\ba{c}Y_{k-1}(:,2:N-k+1)\\A(k+r,k+1:N)\ea\right]
\end{equation}
and determine the $N-k$-row $X_k$ and the $r\times(N-k)$ matrix $Y_k$ from the
partition
\be\label{olz}
Z_k=\left[\ba{c}X_k\\Y_k\end{array}\right].
\end{equation}
Set
\be\label{olia}
R(k,k)=x_k,\quad R(k,k+1:N)=X_k.
\end{equation}

3. For the $r\times r$ matrix $Y_{N-r}$ compute the QR factorization
\be\label{olk}
Y_{N-r}=\hat U_{N-r+1}T_{N-r+1}
\end{equation}
with a unitary $r\times r$ matrix $U_{N-r+1}$ and an upper triangular matrix
$T_{N-r+1}=R(N-r+1:N,N-r+1:N)$ as follows.

3.1. For $k=N-r+1,\dots,N-1$ perform the following.

Set
\be\label{lnk33}
\Delta_k=Y_{k-1}(:,1)
\end{equation}
and determine an $(N-k+1)\times(N-k+1)$ unitary matrix $U_k$ and a number
$x_k$ such that
\be\label{lenk3}
U_k^*\Delta_k=\left[\ba{c}x_k\\0_{(N-k)\times1}\ea\right].
\end{equation}
Compute the $(N-k+1)\times(N-k)$ matrix
\be\label{ol3}
Z_k=U_k^*Y_{k-1}(:,2:N-k+1)
\end{equation}
and determine the $N-k$-row $X_k$ and the $(N-k)\times(N-k)$ matrix $Y_k$ from
the partition
\be\label{olz3}
Z_k=\left[\ba{c}X_k\\Y_k\end{array}\right].
\end{equation}
Set
\be\label{olia3}
R(k,k)=x_k,\quad R(k,k+1:N)=X_k.
\end{equation}

3.2.
Set $x_N=Y_{N-1},\;R(N,N)=x_N$.
\end{theorem}

The proof is performed in a standard way.

Applying Lemma \ref{UU} to the matrix $U^*$ with the matrix $U$ as in
(\ref{un}), (\ref{un1}) we conclude that $U^*$ is a block lower triangular
matrix and obtain the formulas to determine its lower quasiseparable
generators.

\begin{lemma}\label{AA}
The matrix $U^*$ in (\ref{un}), (\ref{un1}) is a unitary lower Green - upper
band of order $r$ matrix with lower Green generators $p_U(i),q_U(i),a_U(i)\;
i=1,\dots,N-r$ and diagonal entries $d_U(k)\;(k=1,\dots,N-r)$ obtained from the
partitions
\be\label{iof5}
U^*_k=\left[\ba{cc}p_U(k)&d_U(k)\\a_U(k)&q_U(k)\end{array}\right],\quad
k=1,\dots,N-r
\end{equation}
and by setting
\be\label{iof5r}
p_U(N-r+1)=\hat U^*_{N-r+1}.
\end{equation}
\end{lemma}

Now we obtain an algorithm to compute lower Green generators
of the inverse of a lower band matrix. We proceed via the formulas $A=UR$ and
next $A^{-1}=R^{-1}U^*$.
\begin{theorem}\label{UF1}
Let $A$ be a lower band of order $r$ matrix.

Then lower Green generators $p(i),q(i),a(i)\;(i=1,\dots,N-r)$,
$P_{N-r+1}$ of the matrix $A^{-1}$ are obtained as follows.

1. Using the algorithm from Theorem \ref{ICF1} compute the unitary matrices
$U_k\;(k=1,\dots,$
$N-1)$ of the orders $r_k=r$ for $r=1,\dots,N-r$ and
$r_k=N-k+1$ for $k=N-r+1,\dots,N-1$,
as well as diagonal entries and subrows
\be\label{tre}
x_k=R(k,k),\;k=1,\dots,N,\quad X_k=R(k,k+1:N),\;k=1,\dots,N-1
\end{equation}
of the lower triangular matrix $R$.
Determine the lower Green
generators $p_U(i),q_U(i)$,
$a_U(i)\;(i=1,\dots,N-r)$ of the unitary
lower Green-upper band matrix $U^*$ via partitions (\ref{iof5}) and the
matrices $p_U(i)),a_U(i)\;(i=N-r+1,\dots,N-1$ of sizes $1\times r_i,
r_i\times r_{i+1}$ from the partitions
\be\label{iof5.}
U^*_i=\left[\ba{c}p_U(i)\\a_U(i)\end{array}\right],\quad
i=N-r+1,\dots,N-1.
\end{equation}

2. Compute the lower Green generators of the matrix $A^{-1}$ as follows.

2.1. Compute the lower Green generator $P_{N-r+1}$ as follows.
Set
\be\label{rn}
p(N)=P_N=\frac1{x_N},
\end{equation}
and for $k=N-1,\dots N-r+1$ compute
\be\label{mishln}
p(k)=\frac1{x_k}(p_U(k)-X_kP_{k+1}a_U(k))
\end{equation}
\be\label{mashln}
P_k=\left(\ba{c}p(k)\\P_{k+1}a_U(k)\ea\right).
\end{equation}

2.2. Compute the lower Green generators $p(k),q(k),a(k)$ as follows.

2.2.1. Set
\be\label{masl}
q(k)=q_U(k),\;a(k)=a_U(k),\;k=1,\dots,N-r.
\end{equation}

2.2.2. For $k=N-r,\dots 1$ compute
\be\label{mishl}
p(k)=\frac1{x_k}(p_U(k)-X_kP_{k+1}a(k))
\end{equation}
\be\label{mashl}
P_k=\left(\ba{c}p(k)\\P_{k+1}a(k)\ea\right).
\end{equation}
\end{theorem}

{\em Proof.} Set
\be\label{lrskk}
R_k=R(1:k,1:k),\;k=1,\dots,N.
\end{equation}
Using the fact that $R$ is an upper triangular matrix we get
\be\label{loras}
R^{-1}(k,k)=\frac1{x_k},\;k=1,\dots,N,\quad
 R^{-1}(k,k+1:N)=-\frac1{x_k}X_kR_{k+1}^{-1},\;k=1,\dots,N-1.
\end{equation}
Hence it follows that
\be\label{lm}
R_k=\left(\ba{cc}\frac1{x_k}&-\frac1{x_k}X_kR_{k+1}^{-1}\\
0&R^{-1}_{k+1}\ea\right),\;k=1,\dots,N-1.
\end{equation}

The $r\times r$ matrix $P_{N-r+1}$ may be determined via $P_{N-r+1}=Y_{N-r}^{-1}$.
Using (\ref{olk}) we have $P_{N-r+1}=R_{N-r+1}^{-1}U^*_{N-r+1}$.
On the Stage 3 of the algorithm from Theorem \ref{ICF1} we obtain the
factorization
$$
\hat U_{N-r+1}=U_{N-r+1}\left(\ba{cc}1&0\\0&U_{N-r+2}\ea\right)\cdots
\left(\ba{cc}I_{r-2}&0\\0&U_{N-1}\ea\right).
$$
Set
$$
\tilde P_N=R_N^{-1}=\frac1{x_N},\quad
\tilde P_t=R_t^{-1}
\left(\ba{cc}I_{t-N-r}&0\\0&U^*_{N-1}\ea\right)\cdots U^*_t,\;t=N-1,\dots,N-r+1.
$$
It is clear that $\tilde P_{N-r+1}=P_{N-r+1}$, we set $P_t=\tilde P_t,\;
t=N-r+2,\dots,N$. We should prove the relations (\ref{rn}), (\ref{mishln}),
(\ref{mashln}). The equality (\ref{rn}) is clear. Assume that for some $k$ with
$N-r+2\le k\le N$ the relations hold. We have
$$
P_{k-1}=R_{k-1}^{-1}\left(\ba{cc}I_{k-N-r+1}&0\\0&U^*_{N-1}\ea\right)\cdots
\left(\ba{cc}1&0\\0&U^*_k\ea\right)U^*_{k-1}
$$
Using (\ref{lm}) and (\ref{iof5.}) we get
\begin{gather*}
P_{k-1}=\left(\ba{cc}\frac1{x_{k-1}}&-\frac1{x_{k-1}}X_{k-1}R_k^{-1}\\
0&R^{-1}_k\ea\right)\\
\left(\left(\ba{cc}I_{k-N-r+1}&0\\0&U^*_{N-1}\ea\right)\cdots
\left(\ba{cc}1&0\\0&U^*_k\ea\right)\right)
\left(\ba{c}p_U(k-1)\\a_U(k-1)\ea\right)
\end{gather*}
and therefore
\begin{gather*}
P_{k-1}=\left(\ba{cc}\frac1{x_{k-1}}&-\frac1{x_{k-1}}X_{k-1}P_k\\0&P_k\ea\right)
\left(\ba{c}p_U(k-1)\\a_U(k-1)\ea\right)=\\
\left(\ba{c}\frac1{x_{k-1}}(p_U(k-1)-X_{k-1}P_ka_U(k-1))\\
P_ka_U(k-1)\ea\right)
\end{gather*}
which completes the proof of (\ref{rn})-(\ref{mashln}).

Now one should justify Stage 2. We apply  Theorem \ref{ULM} with $B=U^*$ and
$S=R^{-1}$. Using the formulas (\ref{mas}) we obtain (\ref{masl}).
Inserting this in (\ref{mish}), (\ref{mash}) and using the
equality (\ref{iof5r}) we get
\be\label{masho}
p(k)=\frac1{x_k}(p_U(k)-X_kR_{k+1}^{-1}P^U_{k+1}a_U(k)),\quad k=N-r,\dots,1
\end{equation}
with
\be\label{mashov}
P^U_{N-r+1}=U_{N-r+1}^*,\quad
P^U_k=\left(\ba{c}p_U(k)\\P^U_{k+1}a_U(k)\ea\right),\; k=N-r,\dots,1.
\end{equation}
Set $P_k=R_k^{-1}P^U_k,\;k=N-r+1,\dots,1$. It remains to show that the
relations (\ref{mishl}), (\ref{mashl}) hold. Indeed using $a(k)=a_U(k)$ we get
$$
P_k=\left(\ba{cc}\frac1{x_k}&-\frac1{x_k}X_kR_{k+1}^{-1}\\
0&R_{k+1}^{-1}\ea\right)\left(\ba{c}p_U(k)\\P^U_{k+1}a(k)\ea\right)=
\left(\ba{c}p(k)\\P_{k+1}a(k)\ea\right),\;k=N-r,\dots,1
$$
with $p(k)$ as in (\ref{mishl}).
$\hfill\Box$

Notice that the lower generators of the matrix $A^{-1}$ obtained in the theorem
are in the right normal form (see \cite[Section 5.8]{EGH1}), i.e.
\be\label{mashen}
a(k)((a(k))^*+q(k)((q(k))^*=I_{r},\;k=2,\dots,N-r.
\end{equation}

\subsection{Two-sided band matrices}

Let us now consider the case of a two-sided band matrix $A$ with half-bandwidth $r$. In other words, $A$
satisfies the conditions $A_{ij}=0$ for $|i-j|>r$. By the Asplund theorem the
inverse matrix $A^{-1}$ is a (two-sided) Green matrix of order $r$, i.e.,
it satisfies both conditions (\ref{natu}), (\ref{natuu}). For such matrices we obtain an inversion algorithm whose arithmetic complexity is linear in $N$.

\begin{theorem}\label{ICF1T}
Let $A=\{A_{ij}\}_{i,j=1}^N$ be a band matrix of order $r$.

The matrix $A$ admits the factorization
\be\label{irn15t}
A=UR,
\end{equation}
where $U$ is a unitary matrix represented as the product (\ref{un}),
(\ref{un1}) involving the $(r+1)\times(r+1)$ unitary matrices
$U_k$, and $R$ is an upper triangular upper band  matrix of order $2r$.
Moreover, the unitary matrices $U_k$ as well as the nonzero upper triangular
entries of the matrix $R$ are obtained as follows.

1.1. Set $\tilde Y_0=A(1:r,1:2r+1)$.

1.2. For $k=1,\dots,N-2r-1$ perform the following.

Determine an $(r+1)\times(r+1)$ unitary matrix $U_k$ and a number $x_k$ via
(\ref{lnk}), (\ref{lenk}).

Compute the $(r+1)\times  2r$ matrix
$$
\tilde Z_k=U_k^*\left[\ba{cc}\tilde Y_{k-1}(:,2:2r)&0_{r\times1}\\
A(k+r,k+1:k+2r-1)&A(k+r,k+2r)\ea\right]
$$
and determine the $2r$-dimensional row $X_k$ and the $r\times 2r$ matrix $Y_k$
from the partition
$$
\tilde Z_k=\left[\ba{c}\tilde X_k\\\tilde Y_k\end{array}\right].
$$
Set $$R(k,k)=x_k,\quad R(k,k+1:k+2r)=\tilde X_k$$.

1.3. For $k=N-2r,\dots,N-r-1$ perform Step 2 of the algorithm from Theorem
\ref{ICF1} to compute $X_k$.

1.4. Perform Step 3 of the algorithm from Theorem \ref{ICF1} to compute
$U_{N-r+1}$,
$Z_{N-r+1}$.

Next, lower Green generators of the matrix $A^{-1}=R^{-1}U^*$ are obtained as
follows.

2.1. Set
\be\label{maslu}
q(k)=q_U(k),\;a(k)=a_U(k),\;k=1,\dots,N-r.
\end{equation}

2.2.1. Set
$$
p(N-r+1)=Z_{N-r+1}^{-1}U^*_{N-r-1},\quad P(N-r+1)=p(N-r+1)
$$
and for $k=N-r,\dots, N-2r+1$ compute
\be\label{mishlu}
p(k)=\frac1{x_k}(p_U(k)-\tilde X_kP_{k+1}a(k))
\end{equation}
\be\label{mashlu}
P_k=\left(\ba{c}p(k)\\P_{k+1}a(k)\ea\right).
\end{equation}

2.2.2. Set
\be\label{msn}
t_{N-2r+1}=P_{N-2r+1}
\end{equation}
and for $k=N-2r,\dots,1$ compute
\be\label{mishluc}
p(k)=\frac1{x_k}(p_U(k)-\tilde X_kt_{k+1}a(k))
\end{equation}
\be\label{mashluc}
t_k=\left(\ba{c}p(k)\\t_{k+1}(1:2r-1,:)a(k)\ea\right).
\end{equation}
\end{theorem}

{\em Proof.}
Since $A$ is an upper band of order $r$ matrix, the matrix $R$ is upper band
of order $2r$. More precisely, inserting
$$
A(k,k+r+1:N)=0,\;k=1,\dots,N-r-1
$$
in (\ref{ol})  we get
$$
X_k=\left(\ba{cc}\tilde X_k&0_{1\times (N-2r-k))}\ea\right),\;k=1,\dots,N-2r-1.
$$
Inserting this in (\ref{mishl}) we obtain (\ref{mishluc}) with
$$
t_k=P_k(1:2r,:)={\rm col}\left(p(j)a^>_{j,k-1}\right)_{j=k}^{k+2r-1},\;k=N-2r,\dots,1.
$$
We show that the values $t_k$ satisfy the recursive relations (\ref{msn}),
(\ref{mashluc}). Indeed (\ref{msn}) is clearly satisfied, and we have
\begin{gather*}
t_{k-1}={\rm col}\left(p(j)a^>_{j,k-2}\right)^{k-2+2r}_{j=k-1}=\\
\left(\ba{c}p(k-1)\\{\rm col}(p(j)a^>_{j,k-1})^{k-2+2r}_{j=k}a(k-1)\ea\right)=
\left(\ba{c}p(k-1)\\t_k(1:2r-1,:)a(k-1)\ea\right)
\end{gather*}
which completes the proof.
$\hfill\Box$

\section{The LU algorithm}\label{sec:LU}
\setcounter{equation}{0}

Let $A$ be a lower band matrix of order $r$. Assume that $A$ is strongly
regular, that is, all its principal minors are nonzero. We derive here an analogue of the results in the previous section, using
the representation of the matrix $A$ in the form $A=LR$ with a unit-lower
triangular matrix $L$ and an upper triangular
matrix $R$. An analogue of Theorem \ref{ICF1} looks as follows.

\begin{theorem}\label{ICF2}
Let $A=\{A_{ij}\}_{i,j=1}^N$ be a strongly regular lower band matrix of order $r$.

Then in the factorization of $A$ as
\be\label{irn15l}
A=LR,
\end{equation}
where $L$ is a unit lower triangular matrix and $R$ is an upper triangular
matrix $R$, the inverse $L^{-1}$ of the lower
triangular factor $L$ may be represented as the product
\be\label{unl}
L^{-1}=\tilde L_{N-r+1}\tilde L_{N-r}\tilde L_{N-r-1}\cdots\tilde L_1
\end{equation}
with
\be\label{un1l}
\tilde L_k=I_{k-1}\oplus L_k\oplus I_{N-k-r},\; k=1,\dots,N-r,\quad
\tilde L_{N-r+1}=I_{N-r}\oplus L_{N-r+1}
\end{equation}
and $(r+1)\times(r+1)$ lower triangular matrices $L_k$ and $r\times r$ lower
triangular matrix $L_{N-r+1}$. Moreover, the lower
triangular matrices $L_k$ as well as the upper triangular
entries of the matrix $R$ are obtained as follows.

1. Set $\gamma_1=A(1,1)$ and compute $f_1=A(2:r+1,1)/\gamma_1$. Set
$X_1=A((1,2:N)$ and compute
$$
Y_1=-f_1\cdot X_1+A(2:r+1,2:N).
$$
Set
$$
L_1=\left(\ba{cc}1&0\\-f_1&I_r\ea\right).
$$
Set $R(1,1)=\gamma_1,\;R(1,2:N)=A(1,2:N)$.

2. For $k=2,\dots,N-r$ perform the following.

Set $\gamma_k=Y_{k-1}(1,1),\;X_k=Y_{k-1}(1,2:N-k+1)$ and compute
$$
f_k=\left(\ba{c}Y_{k-1}(2:r,1)\\A(k+r,k)\ea\right)\frac1{\gamma_k}.
$$
Set
$$
Z_k=\left[\ba{c}Y_{k-1}(2:r,2:N-k+1)\\A(k+r,k+1:N)\ea\right]
$$
and compute
$$
Y_k=-f_k\cdot X_k+Z_k.
$$
Set
$$
L_k=\left(\ba{cc}1&0\\-f_k&I_r\ea\right).
$$
Set $$R(k,k)=\gamma_k,\quad R(k,k+1:N)=X_k.$$


3. Compute the LU factorization
$$
Y_{N-r}=T_{N-r+1}S_{N-r+1},
$$
and set
$$
L_{N-r+1}=T^{-1}_{N-r+1}
$$
with a unit lower triangular matrix $L_{N-r+1}$ ($T_{N-r+1}$) and an upper
triangular matrix $S_{N-r+1}$.
\end{theorem}

This is in fact a structured version of the Gaussian elimination algorithm.

\begin{cor}\label{nats}
In the conditions of Theorem \ref{ICF2}, the unit lower triangular matrix $L$ is
determined elementwise by the formulas
\begin{gather*}
L(i,i+1:i+r)=f_i,\;i=1,\dots,N-r,\quad L(N-r+1:N,N-r+1:N)=T_{N-r+1},\\
 L(i,j)=0,\;i-j>r\;\mbox{or}\;i<j,\quad L(i,i)=1,\;i=1,\dots,N.
\end{gather*}
\end{cor}

\begin{lemma}\label{pl}
Let $L$ be an r-band unit lower triangular matrix as in Theorem \ref{ICF2}
represented elementwise. Set $g_k=L(k+r,k:k+r-1),\;k=N-r,\dots,2$.

Then the matrices $L_k$ in (\ref{unl}), (\ref{un1l}) may be expressed via relations
\be\label{ds}
L_k=\left(\ba{cc}I_r&0\\-g_k&1\ea\right),\;k=1,\dots,N-r,\quad
L_{N-r+1}=((L(N-r+1:N,N-r+1:N))^{-1}.
\end{equation}
\end{lemma}

{\em Proof.} One can check directly that the equality
$$
L\cdot\tilde L_{N-r+1}
\tilde L_{N-r}\tilde L_{N-r-1}\cdots \tilde L_2\tilde L_1=I
$$
with $\tilde L_k$ as in (\ref{un1l}) and $L_k$ as in (\ref{ds}) holds.
$\hfill\Box$

\begin{remark}
In the conditions of Theorem \ref{UF1l}  set
\be\label{mar}
f_k=\left(\ba{c}c_k\\\alpha_k\ea\right),\;
g_k=\left(\ba{cc}\beta_k&b_k\ea\right),\quad k=2,\dots,N-r-1.
\end{equation}
with the columns $c_k$ and rows $b_k$ of the size $r-1$ and the numbers
$\alpha_k,\beta_k$.

Then the equalities $\alpha_k=\beta_k,;k=2,\dots,N-r-1$ hold.
\end{remark}

{\em Proof.} Using Corollary \ref{nats} and the definition of $g_k$ we have
$$
\alpha_k=L(k+r,k)=\beta_k,\;k=2,\dots,N-r-1.\quad
\hfill\Box$$

Applying Lemma \ref{UU} to the matrix $L^{-1}$ we conclude that
$L^{-1}$ is a lower triangular and lower Green
matrix and we obtain the formulas to determine its lower Green
generators.

\begin{lemma}\label{AAL}
The matrix $L^{-1}$ in (\ref{unl}), (\ref{un1l}) is a lower Green - upper band
matrix, whose lower Green generators $p_L(i),q_L(i),a_L(i),\;
i=1,\dots,N-r,\;p_L(N-r+1)$ of orders equal to $r$ and
diagonal entries $d_L(k)\;(k=1,\dots,N-r)$ are obtained from the  partitions
\be\label{iof5l}
L_i=\left[\ba{cc}p_L(i)&d_L(i)\\a_L(i)&q_L(i)\end{array}\right],\quad
i=1,\dots,N-r,
\end{equation}
and by setting
\be\label{iof5ll}
p_L(N-r+1)=L_{N-r+1}.
\end{equation}
\end{lemma}

Now we obtain an algorithm to compute lower Green generators
of the inverse of a lower band strongly regular matrix.

\begin{theorem}\label{UF1l}
Let $A$ be a strongly regular, lower banded matrix of order $r$.

Then lower Green generators $p(i),q(i),a(i),\;i=1,\dots,N-r,\;p(N-r+1)$
of the matrix $A^{-1}$ are obtained as follows.

1. Using the algorithm from Theorem \ref{ICF1} and the formulas from Lemma
\ref{AA} compute the upper triangular matrix $R$ and the lower Green
generators $p_L(i),q_L(i),a_L(i)$,
$i=1,\dots,N-r,\;p_L(N-r+1)$ 
of the lower
triangular matrix $L^{-1}$, as well as diagonal entries and subrows
\be\label{trel}
x_k=R(k,k),\;X_k=R(k,k+1:N),\;k=1,\dots,N-r-1
\end{equation}
and the submatrix
\be\label{trill}
Z_{N-r+1}=R(N-r+1:N,N-r+1:N)
\end{equation}
of the upper triangular matrix $R$ such that $A=LR$.

2. Compute  lower Green generators as follows

2.1. Set
\be\label{masll}
q(k)=q_L(k),\quad a(k)=a_L(k),\quad k=1,\dots,N-r.
\end{equation}

2.2. Set
$$
p(N-r+1)=Z_{N-r+1}^{-1}L_{N-r+1},\quad P_{N-r+1}=p(N-r+1)
$$
and for $k=N-r,\dots, 1$ compute
\be\label{mishll}
p(k)=\frac1{x_k}(p_L(k)-X_kP_{k+1}a(k))
\end{equation}
\be\label{mashll}
P_k=\left(\ba{c}p(k)\\P_{k+1}a(k)\ea\right).
\end{equation}
\end{theorem}

\begin{lemma}
In the conditions of Theorem \ref{UF1l}  the lower Green  generators
$a(k),q(k),$
$k=2,\dots,N-r$ are given by the formulas
$$
a(k)=\left(\ba{cc}-c_k&I_{r-1}\\-\alpha_k&0\ea\right),\;
q(k)=\left(\ba{c}0_{(r-1)\times1}\\1\ea\right),\quad k=1,\dots,N-r
$$
or by the formulas
$$
a(k)=\left(\ba{cc}0&I_{r-1}\\-\alpha_k&-b_k\ea\right),\;
q(k)=\left(\ba{c}0_{(r-1)\times1}\\1\ea\right),\quad k=1,\dots,N-r,
$$
with $c_k,b_k,\alpha_k$ as in (\ref{mar}).
\end{lemma}

The analogue of Theorem \ref{ICF1T} looks as follows.

\begin{theorem}\label{I1CF1L}
Let $A=\{A_{ij}\}_{i,j=1}^N$ be a strongly regular, two-sided band matrix of order $r$.

The matrix $A$ admits the factorization
\be\label{irn15ll}
A=LR,
\end{equation}
where $L$ is a lower triangular matrix represented as the product (\ref{un}),
(\ref{un1}) with the lower triangular matrices
$L_k$, and $R$ is an upper triangular upper band  matrix of order $r$.
Moreover, the lower triangular matrices $L_k$ as well as the nonzero upper
triangular entries of the matrix $R$ are obtained as follows.

1.1. Set $\gamma_1=A(1,1)$ and compute $f_1=A(2:r+1,1)/\gamma_1$. Set
$X_1=A((1,2:r+1)$ and compute
$$
Y_1=-f_1\cdot\left[\ba{cc}X_1&0\ea\right]+A(2:r+1,2:r+2).
$$
Set
$$
L_1=\left(\ba{cc}1&0\\-f_1&I_r\ea\right).
$$
Set $R(1,1)=\gamma_1,\;R(1,2:r+1)=A(1,2:r+1)$.

1.2. For $k=2,\dots,N-r-1$ perform the following.

Set $\gamma_k=Y_{k-1}(1,1),\;X_k=Y_{k-1}(1,2:r+1)$ and compute
$$
f_k=\left(\ba{c}Y_{k-1}(2:r,1)\\A(k+r,k)\ea\right)\frac1{\gamma_k}.
$$
Set
$$
Z_k=\left[\ba{c}Y_{k-1}(2:r,2:r+1)\,\,\, A(k+1:k+r-1,k+r+1)\\A(k+r,k+1:k+r+1)\ea\right]
$$
and compute
$$
Y_k=-f_k\cdot\left[\ba{cc}X_k&0\ea\right]+Z_k.
$$
Set
$$
L_k=\left(\ba{cc}1&0\\-f_k&I_r\ea\right).
$$
Set $$R(k,k)=\gamma_k,\quad R(k,k+1:k+r)=X_k.$$

1.3.
Set $\gamma_{N-r}=Y_{N-r-1}(1,1),\;X_{N-r}=Y_{N-r-1}(1,2:r+1)$ and compute
$$
f_{N-r}=\left(\ba{c}Y_{N-r-1}(2:r,1)\\A(N,N-r)\ea\right)\frac1{\gamma_{N-r}}.
$$
Set
$$
Z_{N-r}=\left[\ba{c}Y_{N-r-1}(2:r,2:r+1)\\A(N,N-r+1:N)\ea\right]
$$
and compute
$$
Y_{N-r}=-f_{N-r}X_{N-r}+Z_{N-r}.
$$
Set
$$
L_{N-r}=\left(\ba{cc}1&0\\-f_{N-r}&I_r\ea\right).
$$
Set $$R(N-r,N-r)=\gamma_{N-r},\quad R(N-r,N-r+1:N)=X_{N-r}.$$


Compute the LU factorization
$$
Y_{N-r}=T_{N-r+1}S_{N-r+1},
$$
and set
$$
L_{N-r+1}=T^{-1}_{N-r+1}, \qquad R(N-r+1:N,N-r+1:N) = S_{N-r+1}
$$
with a unit lower triangular matrix $L_{N-r+1}$ ($T_{N-r+1}$) and an upper
triangular matrix $S_{N-r+1}$.

Next, lower Green generators of the matrix $A^{-1}=R^{-1}L$ are obtained as
follows.

2.1. Set
\be\label{maslul}
q(k)=q_L(k),\;a(k)=a_L(k),\;k=1,\dots,N-r.
\end{equation}

2.2. Set
$$
p(N-r+1)=S_{N-r+1}^{-1}L_{N-r+1},\; P_{N-r+1}=p(N-r+1),\quad
t_{N-r+1}=p(N-r+1)
$$
and for $k=N-r,\dots,1$ compute
\be\label{mishlucl}
p(k)=\frac1{x_k}\left(p_L(k)- X_kt_{k+1}a(k)\right)
\end{equation}
\be\label{mashlucl}
t_k=\left(\ba{c}p(k)\\t_{k+1}(1:r-1,:)a(k)\ea\right).
\end{equation}
\end{theorem}

\section{Numerical tests}\label{sec:test}

We propose a few numerical test that highlight the computational properties of the inversion algorithms presented in Theorems \ref{UF1}, \ref{ICF1T}, \ref{UF1l}  and \ref{I1CF1L}.\footnote{The MATLAB implementation used in these tests is available from {\tt https://people.cs.dm.unipi.it/boito/Green.zip}}

{\bf Example 1.} Here $A$ is an $N\times N$ banded matrix with bandwidth $r=5$ and  random entries, defined via the MATLAB command {\tt A = tril(triu(rand(N),-r),r);}

\begin{figure}
\begin{center}
\caption{Logarithmic plot of running times for Example 1 (red dots), i.e., QR-based inversion of two-sided banded matrices. Matrix size ranges from $250$ to $2000$. The slope of the linear fit is $0.9373$, therefore consistent with theoretical complexity analysis.}\label{fig:QRtimes}
\includegraphics[width=0.8\textwidth]{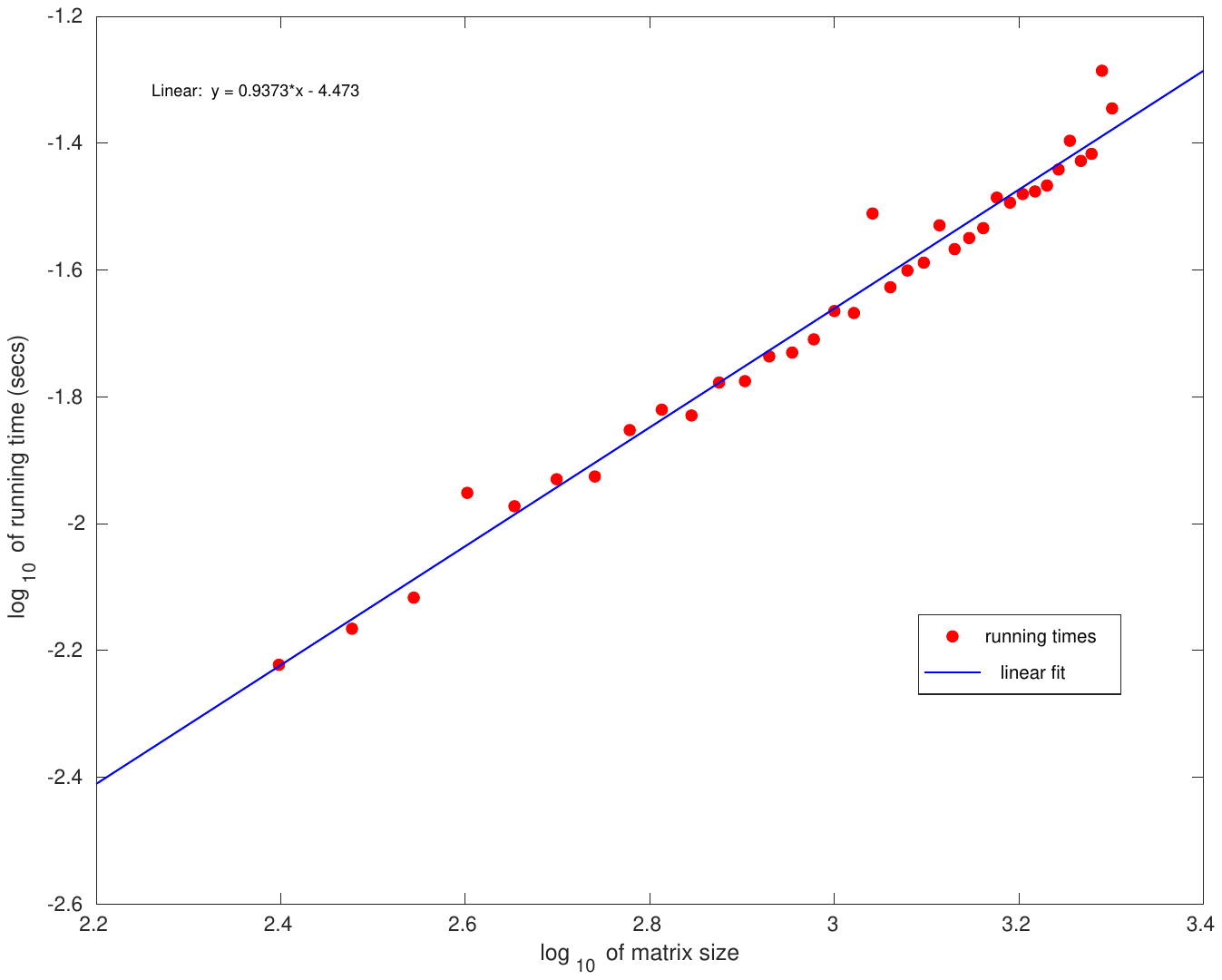}
\end{center}
\end{figure}

\begin{figure}
\begin{center}
\caption{Relative forward errors for Example 1 (red dots), i.e., QR-based inversion of two-sided banded matrices. Matrix size ranges from $250$ to $2000$. Experimental accuracy is consistent with theoretical estimates (black stars) given by the machine epsilon times the 2-norm condition number, and does not appear to deteriorate with increasing size.}\label{fig:QRerr}
\includegraphics[width=0.8\textwidth]{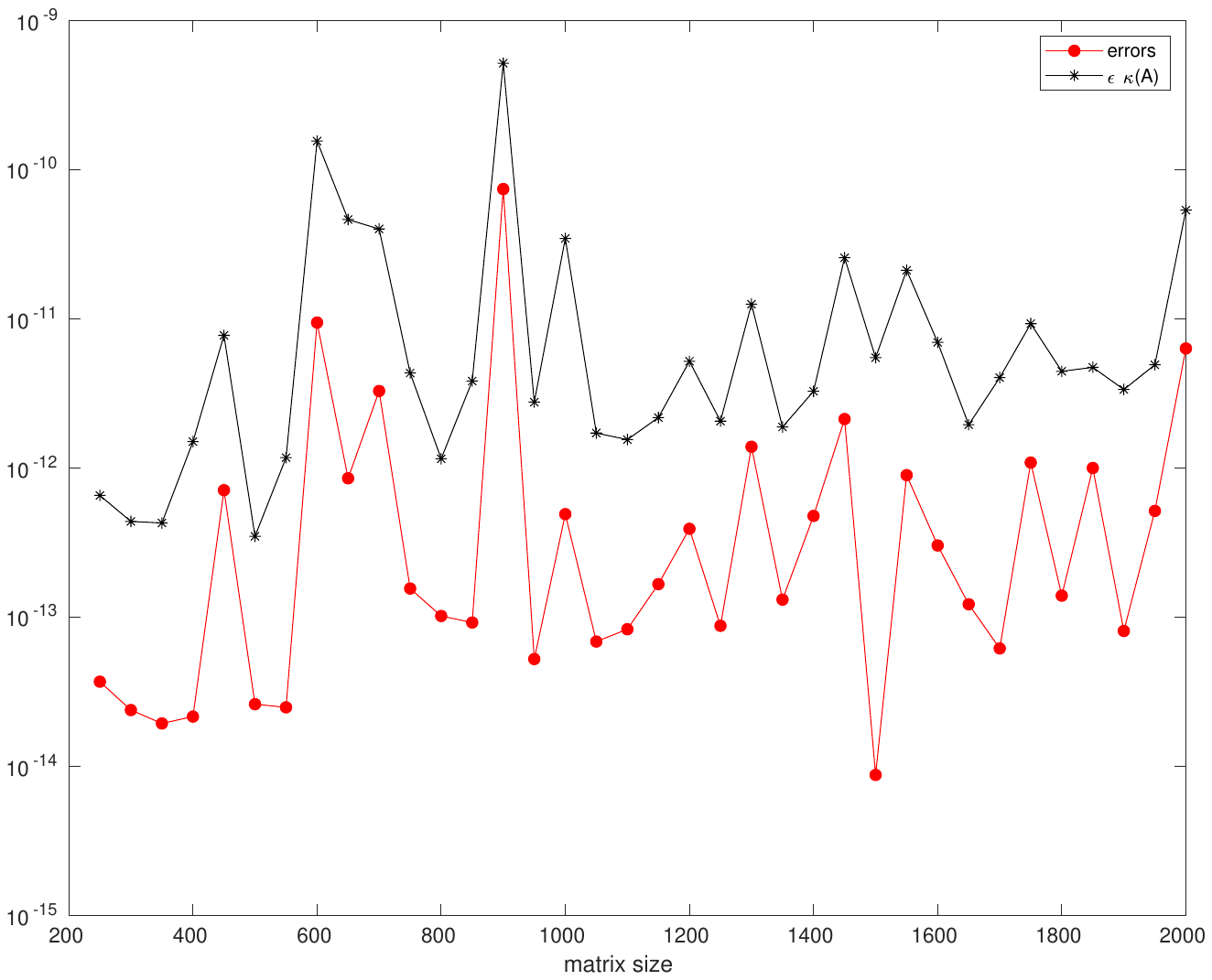}
\end{center}
\end{figure}

Figures \ref{fig:QRtimes} and \ref{fig:QRerr} show running times and relative errors when the QR-based inversion algorithm from Theorem \ref{ICF1T} is applied to $A$, for $N$ ranging between $250$ and $2000$. Here and in the next examples, relative errors are computed w.r.t.~the upper-banded, rank-structured portion of $A^{-1}$ parameterized by lower Green generators, that is, via the MATLAB command
\begin{center}
{\tt err = norm(tril(B,r-1)-tril(Ainv,r-1))/norm(tril(Ainv,r-1))}
\end{center}
where {\tt Ainv} is the inverse of $A$ computed using the MATLAB command {	\tt inv}, and {\tt B} is the portion of $A^{-1}$ reconstructed from lower Green generators provided by Theorem \ref{ICF1T}.

A linear fit on the log-log time plot confirms that the algorithm has arithmetic complexity $O(N)$. The relative forward errors are consistent with the usual estimate $\varepsilon\kappa_2(A)$.

{\bf Example 2.} In this example we test the LU-based inversion algorithm from Theorem \ref{I1CF1L}. The setup is very similar to Example 1, except that here we add a diagonal term to $A$ to ensure strong regularity: we set
\begin{center}
{\tt A = tril(triu(rand(N),-r),r) + r*eye(N);}
\end{center}
again with $r=5$ and $N$ ranging between $500$ and $2500$.

\begin{figure}
\begin{center}
\caption{Logarithmic plot of running times for Example 2 (blue dots), i.e., LU-based inversion of two-sided banded matrices. Matrix size ranges from $500$ to $2500$. The slope of the linear fit is $1.016$, therefore consistent with theoretical complexity analysis.}\label{fig:LUtimes}
\includegraphics[width=0.8\textwidth]{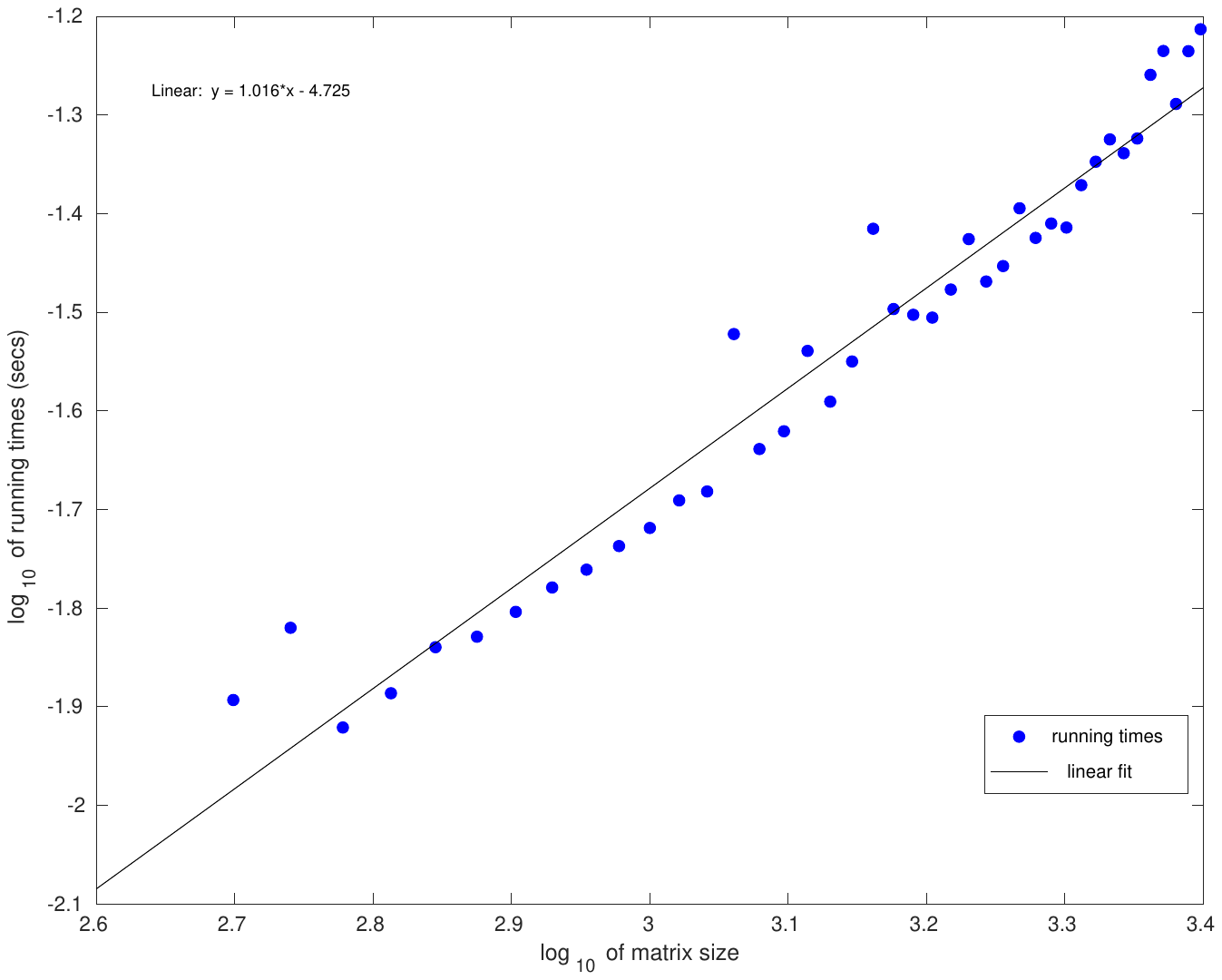}
\end{center}
\end{figure}

\begin{figure}
\begin{center}
\caption{Relative forward errors for Example 2, i.e., LU-based inversion of two-sided banded matrices. Matrix size ranges from $500$ to $2500$.}\label{fig:LUerr}
\includegraphics[width=0.8\textwidth]{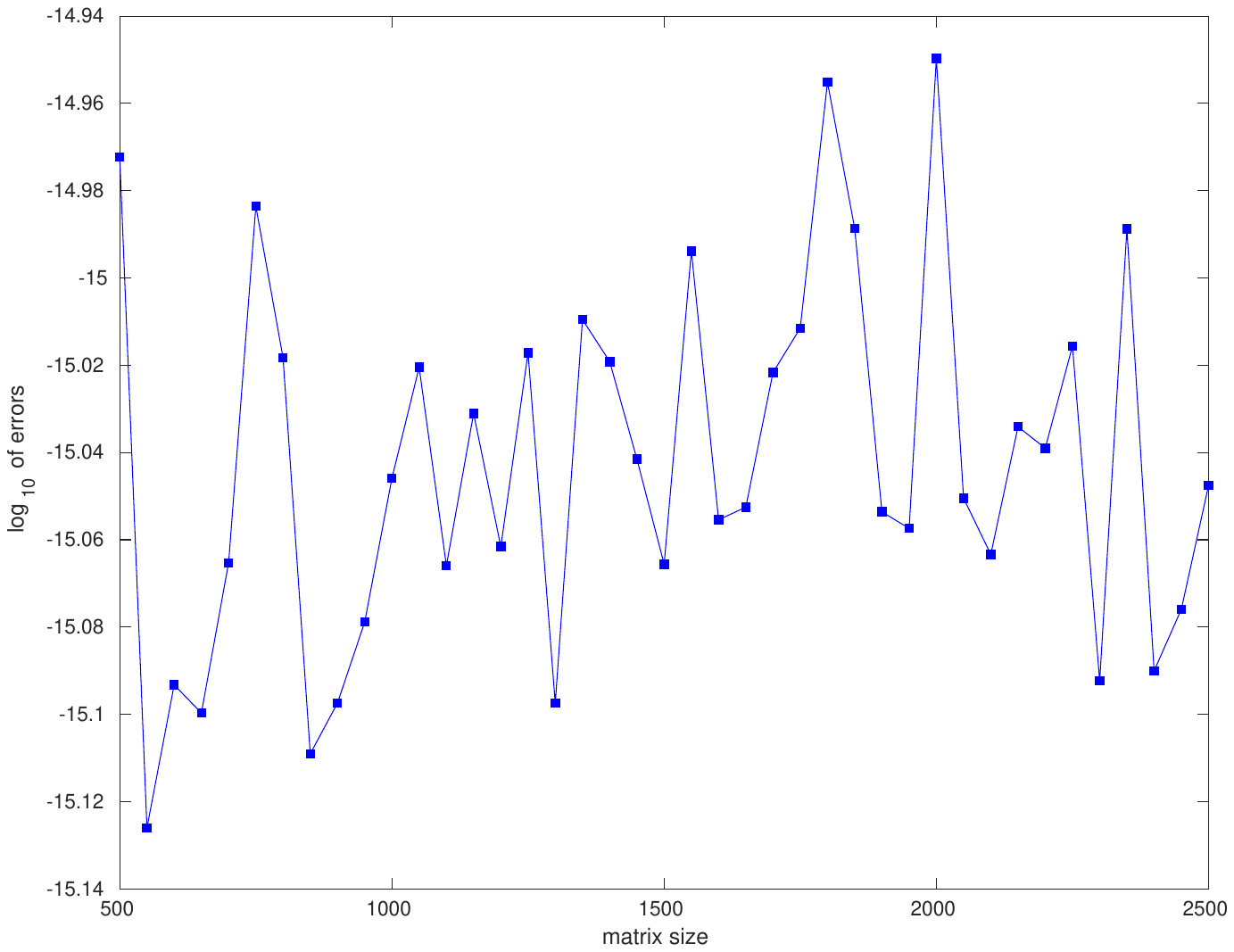}
\end{center}
\end{figure}

Figures \ref{fig:LUtimes} and \ref{fig:LUerr} show running times and relative errors, together with a linear fit on the log-log time plot. The relative forward errors are consistently small. Note that the matrices used in this test are all well-conditioned: their $2$-norm condition number does not exceed $4$.

\begin{figure}
\begin{center}
\caption{Logarithmic plot of running times for Example 3 (two-sided banded case). Matrix size ranges from $600$ to $2500$. Linear fits (not shown here) have a slope of about $1.04$ for both quasiseparable algorithms, $2.40$ for {\tt sparseinv}, $2.75$ for {\tt inv} and $2.78$ for classical QR.}\label{fig:timec2}
\includegraphics[width=0.8\textwidth]{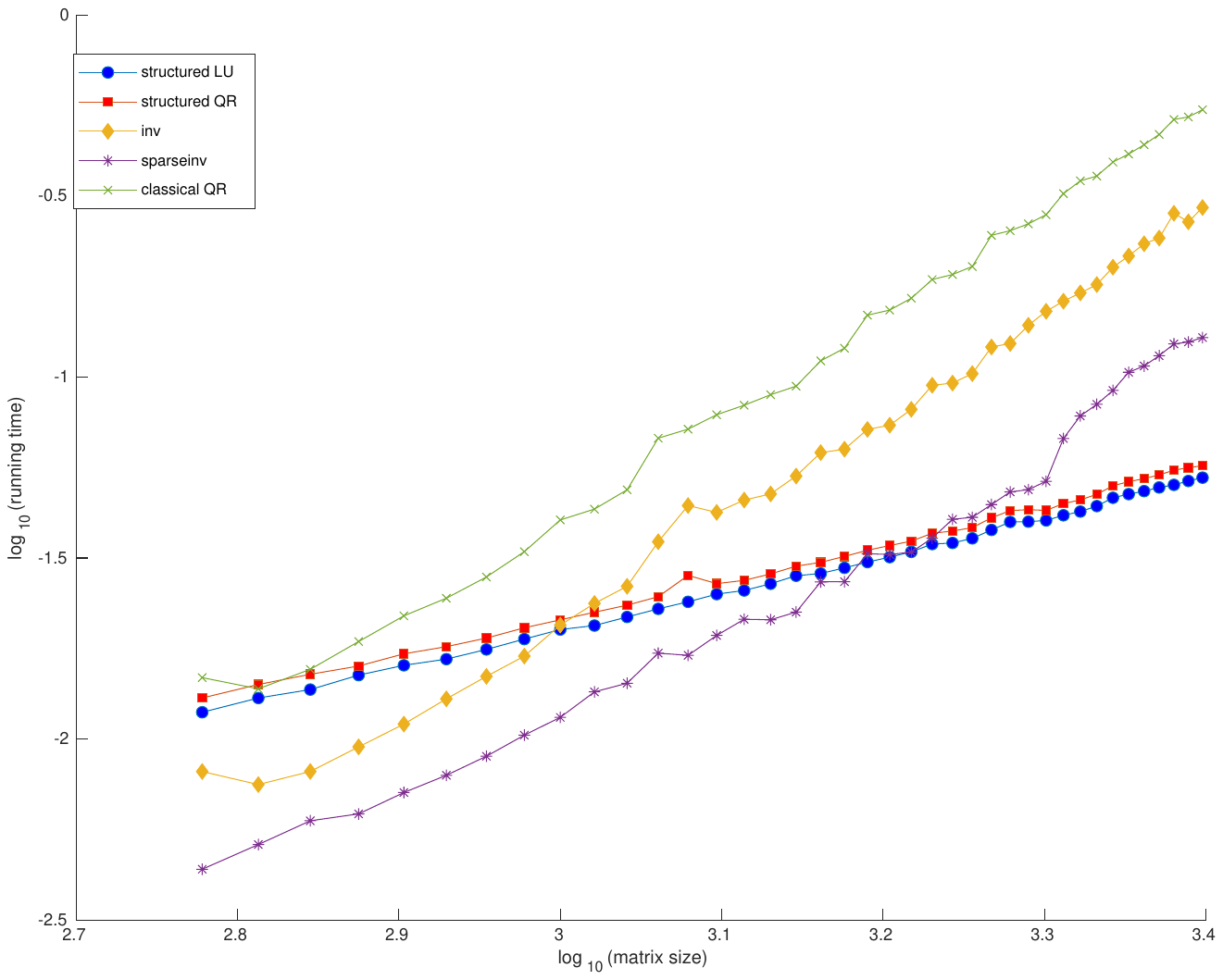}
\end{center}
\end{figure}

\begin{figure}
\begin{center}
\caption{Logarithmic plot of running times for Example 3 (lower banded case). Matrix size ranges from $500$ to $2500$. Linear fits (not shown here) have a slope of about $1.73$ for structured LU, $1.83$ for structured QR, $3.36$ for {\tt sparseinv}, $2.64$ for {\tt inv} and $3.10$ for classical QR.}\label{fig:timec1}
\includegraphics[width=0.8\textwidth]{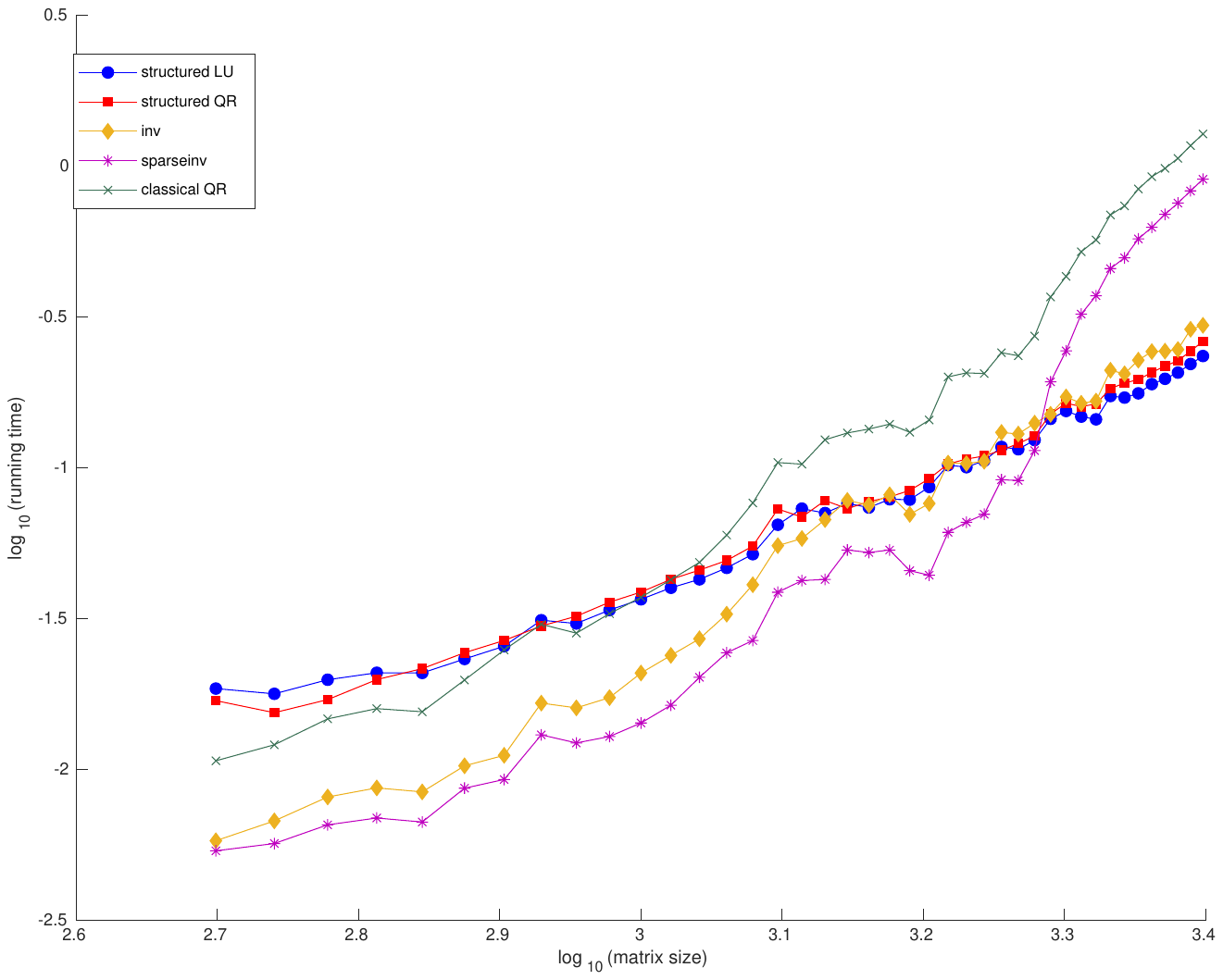}
\end{center}
\end{figure}

{\bf Example 3.} Let us also compare running times for the structured algorithms presented here and for their classical counterparts, namely, MATLAB's command {\tt inv} (which is based on LU factorization), its sparse version {\tt sparseinv} and classical QR-based inversion (using MATLAB's {\tt qr}). Figure \ref{fig:timec1} shows results for lower banded matrices, whereas timings for two-sided banded matrices are plotted in Figure \ref{fig:timec2}. Note that we are comparing here a MATLAB implementation with built-in functions: such a setup clearly penalizes the structured methods. On the other hand, one could argue that the classical methods provide all the entries of the matrix inverse, whereas the structured algorithms  compute quasiseparable generators; if one wanted to reconstruct the whole matrix, clearly the linear complexity would be lost. However, quasiseparable generators are well-suited to many applications, e.g., cases where one needs only a few matrix entries, or where quasiseparable generators are actually needed to carry on further structured computations. The development of decay bounds mentioned in the Introduction is another example.

\begin{figure}
\begin{center}
\caption{Relative errors for Example 4: red dots are for QR-based inversion, blue squares for the classical algorithm provided by the MATLAB command {\tt inv}. Black stars denote theoretical estimates.}\label{fig:QRcond}
\includegraphics[width=0.8\textwidth]{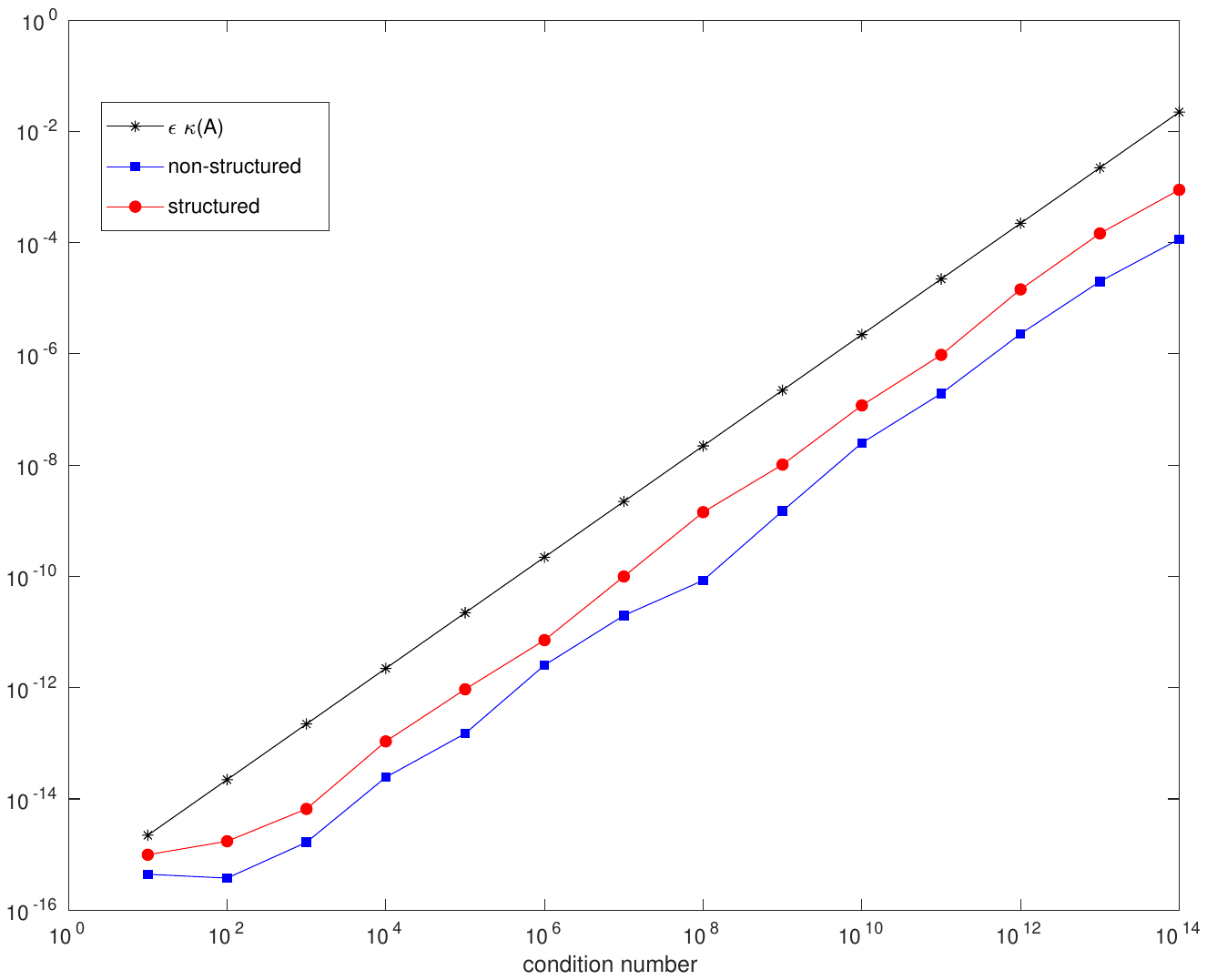}
\end{center}
\end{figure}

{\bf Example 4.} This is a test for the QR-based inversion algorithm for lower-banded matrices proposed in Theorem \ref{UF1}. We focus on the behavior of the algorithm for ill-conditioned matrices. To this end, we form fourteen $100\times 100$ random matrices with lower bandwidth $r=5$ and prescribed condition numbers $10^c$, for $c=1,\ldots,14$. Figure \ref{fig:QRcond} shows relative errors on the rank-structured part of $A^{-1}$, both for the structured algorithm of Theorem \ref{UF1} and for the MATLAB command {\tt inv}. Errors are defined here w.r.t.~the inverse of $A$ computed via variable precision arithmetic (VPA) with 64 digits. Since we are dealing with randomly generated matrices, for each choice of condition number we have performed an average over $10$ cases.

The results for the structured algorithm are consistent with theoretical estimates and comparable to the standard algorithm.

\begin{figure}
\begin{center}
\caption{Relative errors for Example 5: blue dots are for LU-based inversion, red squares for QR-based inversion.}\label{fig:unstableLU}
\includegraphics[width=0.8\textwidth]{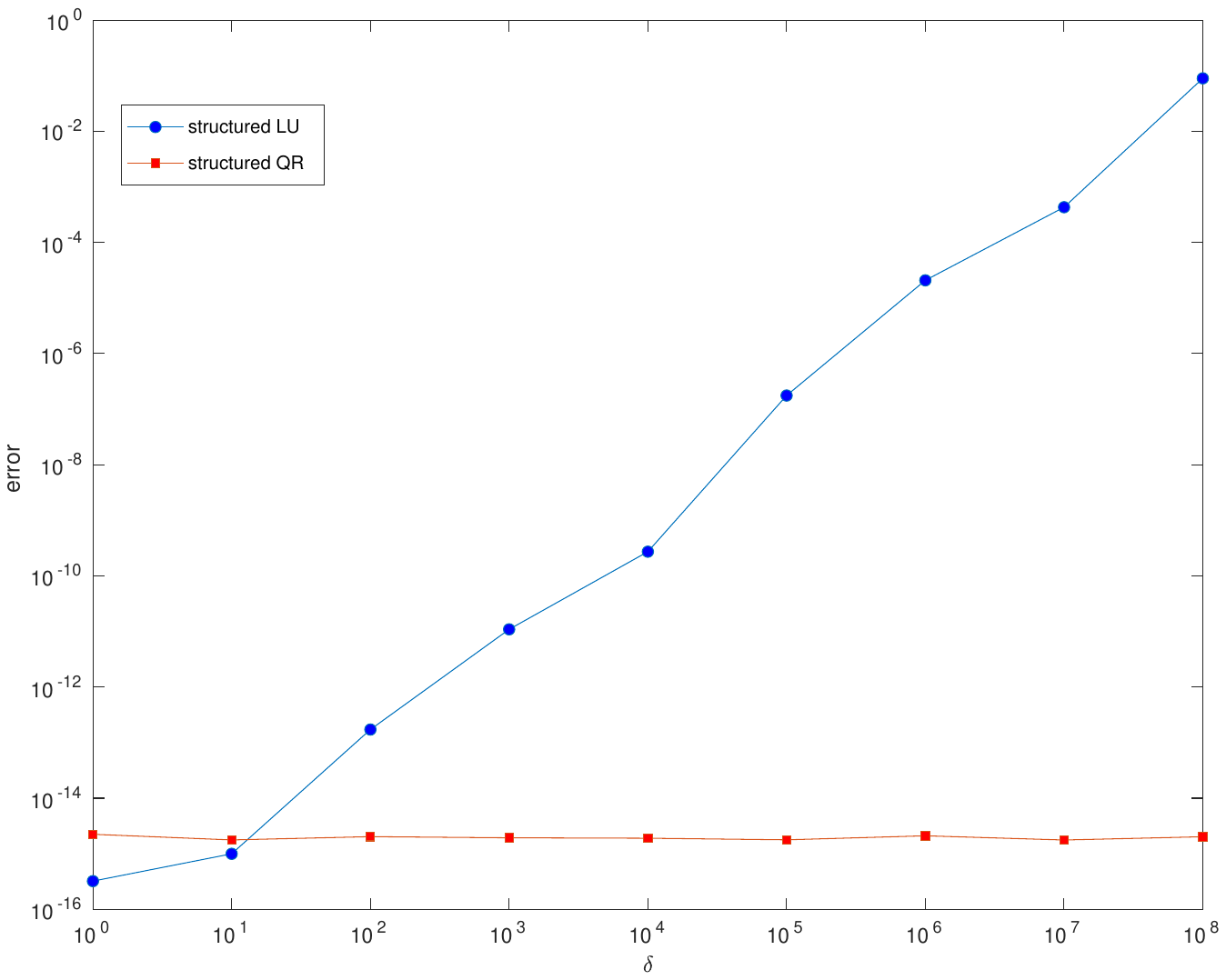}
\end{center}
\end{figure}

{\bf Example 5.} As mentioned above, Theorem \ref{UF1l} is a quasiseparable version of the classical Gauss/LU factorization algorithm. If, for instance, we set $N=10$, $r=2$ and
\begin{center}
{\tt A = triu(rand(N),-r)+r*eye(N);}
\end{center}
we typically obtain a matrix for which the Gauss algorithm is stable. Analogously, the algorithm from Theorem \ref{UF1l} will yield an error of the order of the machine epsilon. However, if we modify the matrix so that small pivots emerge in the Gauss computation, we expect the structured algorithm to exhibit instability, just like the classical algorithm would. In this example we replace the principal $3\times 3$ block of $A$ with
$$
\left[\begin{array}{ccc}
1 & 1 & 1\\
2 & 2+\delta & 5\\
4 & 6 & 8
\end{array}\right],
$$
where $\delta=10^0, 10^{-1},\ldots, 10^{-8}$. The error grows accordingly, as shown in Figure \ref{fig:unstableLU}. On the other hand, if we apply the QR-based inversion algorithm, the error stays small, in accordance with usual error analysis for LU and QR factorizations. Note that $A$ is well-conditioned for all the chosen values of $\delta$.

\section{Conclusions and future work}
In this paper we have proposed inversion algorithms for one- and two-sided banded matrices that rely on a novel quasiseparable formulation of Asplund's theorem. In the two-sided band case, complexity is linear w.r.t.~matrix size, as expected from theory.

Building on explicit inversion formulas presented in Theorems \ref{UF1},  \ref{ICF1T}, \ref{UF1l}, \ref{I1CF1L}, we plan to derive computable {\em a priori} bounds for the off-diagonal decay of inverses of band matrices. Such bounds would rely on a quasiseparable generator representation of band matrices. It will be interesting to compare them to existing bounds based on polynomial approximation of the inverse function as in \cite{DMS84}, which are essentially based on spectral properties of the matrix to be inverted and, because of their generality, can sometimes be quite pessimistic. We expect that quasiseparable-based bounds can be better tailored to the specific matrix under study and can help give a more accurate description of decay behavior.


\begin{thebibliography}{88}

\bibitem{Asplund59}
E.~Asplund, Inverses of matrices $\{a_{ij}\}$ which satisfy $a_{ij}= 0$ for $j> i+ p$. Mathematica Scandinavica (1959), 57--60.

\bibitem{Bapat07}
R.~Bapat,  On generalized inverses of banded matrices. The Electronic Journal of Linear Algebra 16 (2007), 284--290.

\bibitem{BF81}
W.~W.~Barrett and P.~J.~Feinsilver, Inverses of banded matrices. Linear Algebra and its Applications 41 (1981), 111--130.

\bibitem{B91}
R.~Bevilacqua, Structural and computational properties of band matrices, in Complexity of Structured Computational Problems, R.~Bevilacqua, D.~Bini, M.~Capovani, G.~Capriz, B.~Codenotti, M.~Leoncini, G.~Resta, and P.~Zellini, eds., Appl. Math. Monographs, Consiglio Nazionale delle Ricerche, Giardini Editori e Stampatori in Pisa, 1991, pp. 131–188.

\bibitem{BLR90}
R.~Bevilacqua, C.~Lotti, and F.~Romani, Storage compression of inverses of band matrices. Computers and Mathematics with Applications 20 (1990), 1--11.

\bibitem{BenziGolub99}
M.~Benzi and G.~H.~Golub, Bounds for the entries of matrix functions with applications to preconditioning. BIT Numerical Mathematics 39 (1999), 417--438.

\bibitem{BenziRazouk07}
M.~Benzi, M. and N.~Razouk, Decay bounds and $O(n)$ algorithms for approximating functions of sparse matrices. Electron. Trans. Numer. Anal 28 (2007), 16--39).

\bibitem{BenziRinelli}
M.~Benzi and M.~Rinelli, Refined decay bounds on the entries of spectral projectors associated with sparse Hermitian matrices. Linear Algebra and its Applications 647 (2022), 1--30.

\bibitem{Bevilacqua05}
R.~Bevilacqua, E.~Bozzo, G.~M.~ Del Corso, and D.~Fasino, Rank structure of generalized inverses of rectangular banded matrices. Calcolo 42 (2005), 157--169.

\bibitem{BiniMeini99}
D.~A.~Bini and B.~Meini, Effective methods for solving banded Toeplitz systems. SIAM Journal on Matrix Analysis and Applications 20 (1999), 700--719.

\bibitem{BuenoFurtado}
M.~I.~Bueno and S.~Furtado, Singular matrices whose Moore-Penrose inverse is tridiagonal. Applied Mathematics and Computation (2023): 128154.

\bibitem{Canuto14}
C.~Canuto, V.~Simoncini, and M.~Verani, On the decay of the inverse of matrices that are sum of Kronecker products. Linear Algebra and its Applications 452 (2014), 21--39.

\bibitem{Capovani70}
M.~Capovani, Sulla determinazione della inversa delle matrici tridiagonali e tridiagonali a blocchi. Calcolo 7 (1970), 295–303

\bibitem{CG03}
S.~Chandrasekaran and M.~Gu, Fast and stable algorithms for banded plus semiseparable
systems of linear equations. SIAM J. Matrix Anal. Appl., 25 (2003), 373–384

\bibitem{DMS84}
S.~Demko, W.~F.Moss, and P.~W.~Smith,
Decay rates for inverses of band matrices. Mathematics of Computation 43.168 (1984), 491-499.

\bibitem{EGH1}
Y.~Eidelman, I.~Gohberg, and I.~Haimovici,
Separable type representations of matrices and fast algorithms.
Volume 1. Basics. Completion problems. Multiplication and inversion
algorithms, Operator Theory: Advances and Applications,  Birkh\"auser, 2013.

\bibitem{Fasino02}
D.~Fasino and L.~Gemignani, Structural and computational properties of possibly singular semiseparable matrices. Linear Algebra and its Applications 340 (2002), 183--198.

\bibitem{Frommer18}
A.~Frommer, C.~Schimmel, and M.~Schweitzer, Bounds for the decay of the
entries in inverses and Cauchy-Stieltjes functions of certain sparse, normal
matrices. Numer. Linear Algebra Appl. 25 (2018):e2131.

\bibitem{GKbook}
F.~R.~Gantmacher and M.~G.~Krein, Oscillation matrices and kernels and small oscillations of
mechanical systems. (Russian) GITTL, Moscow, 1941. (English translation: AMS, Providence,
2002.)

\bibitem{Garcia2012}
P.~Garc\`\i\mbox{}a-Risue\~no and P.~Echenique, Linearly scaling direct method for accurately inverting sparse banded matrices. Journal of Physics A: Mathematical and Theoretical 45 (2012), 065204.

\bibitem{KS13}
E.~K\i l\i \c{c} and P.~Stanica, The inverse of banded matrices. Journal of Computational and Applied Mathematics 237 (2013), 126--135.

\bibitem{Meurant92}
G.~Meurant, A review on the inverse of symmetric tridiagonal and block tridiagonal matrices. SIAM Journal on Matrix Analysis and Applications 13 (1992), 707--728.

\bibitem{Nabben99}
R.~Nabben, Decay rates of the inverse of nonsymmetric tridiagonal and band matrices. SIAM Journal on Matrix Analysis and Applications 20 (1999), 820--837.

\bibitem{Noschese13}
S.~Noschese, L.~Pasquini, and L.~Reichel, Tridiagonal Toeplitz matrices: properties and novel applications. Numerical Linear Algebra with Applications 20 (2013), 302--326.

\bibitem{Olshevsky10}
V.~Olshevsky, G.~Strang, and P.~Zhlobich, Green’s matrices. Linear Algebra and Its Applications 432 (2010), 218--241.

\bibitem{Romani86}
F.~Romani, On the additive structure of the inverses of banded matrices. Linear Algebra and Its Applications 80 (1986), 131--140.

\bibitem{Rosza91}
P.~R\'ozsa, R.~Bevilacqua, F.~Romani, and P.~Favati. On band matrices and their inverses. Linear Algebra and Its Applications 150 (1991), 287--295.

\bibitem{SSbiblio}
R.~Vandebril, M.~Van Barel, G.~Golub, and N.~Mastronardi. A bibliography on semiseparable matrices. Calcolo 42 (2005), 249--270.

\bibitem{VandebrilBook}
R.~Vandebril, M.~Van Barel, and N.~Mastronardi. Matrix computations and semiseparable
matrices. Vol. 1. Johns Hopkins University Press, Baltimore, MD, 2008. Linear systems.

\bibitem{YI79}
T.~Yamamoto and Y.~Ikebe, Inversion of band matrices. Linear Algebra and Its Applications 24 (1979), 105-111.

\end{thebibliography}
\end{document}